\pdfoutput=1
\documentclass[11pt, a4paper]{article}

\usepackage{fullpage, comment}
\usepackage{latexsym}
\usepackage{import, pst-plot}
\usepackage{amsfonts}
\usepackage{amsmath}
\usepackage{amssymb}
\usepackage{graphicx}

\newcommand{\beq}{\begin{equation}}
\newcommand{\eeq}{\end{equation}}
\newcommand{\ba}{\begin{array}}
\newcommand{\ea}{\end{array}}
\newcommand{\bea}{\begin{eqnarray}}
\newcommand{\eea}{\end{eqnarray}}
\newcommand{\bc}{\begin{center}}
\newcommand{\ec}{\end{center}}
\newcommand{\bt}{\begin{table}}
\newcommand{\et}{\end{table}}

\newcommand{\la}[1]{\label{#1}}
\newcommand{\p}{\partial}

\newcommand{\ds}{\displaystyle}
\newcommand{\no}{\noindent}

\newcommand{\rf}[1]{(\ref{#1})}
\newcommand{\beqno}{\begin{displaymath}}
\newcommand{\eeqno}{\end{displaymath}}

\newcommand{\been}{\begin{enumerate}}
\newcommand{\een}{\end{enumerate}}

\newcommand{\vs}{\vspace*{0.1in}}
\newcommand{\ra}{\rightarrow}
\newcommand{\Ra}{\Rightarrow}
\newcommand{\sgn}{\mbox{sign}}
\newcommand{\twovector}[2]{\left(\ba{c}#1\\#2\ea\right)}
\newcommand{\twomatrix}[4]{\left(\ba{cc}#1 & #2\\#3 & #4\ea\right)}
\renewcommand{\Im}{\mbox{Im}}
\renewcommand{\Re}{\mbox{Re}}

\def\XXint#1#2#3{{\setbox0=\hbox{$#1{#2#3}{\int}$}
     \vcenter{\hbox{$#2#3$}}\kern-.5\wd0}}

\newlength{\myheight}
\newlength{\mylength}

\newcounter{saveeqn}

\newcommand{\alpheqn}{\setcounter{saveeqn}{\value{equation}}
\stepcounter{saveeqn}\setcounter{equation}{0}
\renewcommand{\theequation}{\mbox{\arabic{saveeqn}\alph{equation}}}}

\newcommand{\resetalpheqn}{\setcounter{equation}{\value{saveeqn}}
\renewcommand{\theequation}{\arabic{equation}}}

\newtheorem{example}{Example}

\begin{document}

\title{Fokas's Unified Transform Method for linear systems}

\author{
Bernard Deconinck$^\dagger$, Qi Guo$^*$, Eli Shlizerman$^\#$, and Vishal Vasan$^\ddagger$\\
~\\
$^\dagger$ Department of Applied Mathematics\\
University of Washington\\
Campus Box 353925, Seattle, WA, 98195, USA\\
e-mail: {\tt deconinc@uw.edu, shlizee@uw.edu}\\
~\\
$^*$ UCLA Mathematics Department\\
Box 951555\\
Los Angeles, CA 90095-1555\\
e-mail: {\tt qiguo@math.ucla.edu}\\
~\\
$^\#$ Department of Applied Mathematics \& Department of Electrical Engineering\\
University of Washington\\
Campus Box 353925, Seattle, WA, 98195, USA\\
e-mail: {\tt deconinc@uw.edu, shlizee@uw.edu}\\
~\\
$^\ddagger$ International Centre for Theoretical Sciences\\
Bengaluru, India\\
e-mail: {\tt vishal.vasan@icts.res.in}}

\maketitle

\begin{abstract}
We demonstrate the use of the Unified Transform Method or Method of Fokas for boundary value problems for systems of constant-coefficient linear partial differential equations. We discuss how the apparent branch singularities typically appearing in the global relation are removable, allowing the method to proceed, in essence, as for scalar problems. We illustrate the use of the method with boundary value problems for the Klein-Gordon equation and the linearized Fitzhugh-Nagumo system. The case of wave equations is treated separately in an appendix.
\end{abstract}

\section{Introduction}

The Unified Transform Method (UTM) or Method of Fokas presents a new approach to the solution of boundary value problems (BVPs) for integrable nonlinear equations \cite{fokas1997, fokasbook}, or even more successfully for linear, constant coefficients partial differential equations (PDEs) \cite{dtv, fokasbook}. In this latter context, the method has resulted in new results for one-dimensional problems involving more than two spatial derivatives \cite{fokasima, fokasbook}, for elliptic problems \cite{ashtonfokas, fokasbook}, and for interface problems \cite{deconinckpellonisheils, greeks}, to name but a few. A more complete picture can be obtained from the online repository at \cite{portal}.

Fewer studies exist of how the method applies to linear systems of equations or to higher-order scalar evolution equations that may be rewritten as such. In fact, to our knowledge, the application of the UTM to the standard wave equation presented in the appendix here is new, albeit not surprising. Systems of equations are examined in \cite{fokaspellonisys, pp3, pp1, pp2} and \cite{treharnefokas}. In \cite{treharnefokas}, a problem from thermoelastic deformation is examined, while in \cite{fokaspellonisys} two water wave-related systems are treated. In \cite{pp1}, the UTM is applied to the wave equation in a moving domain, while in \cite{pp3, pp2} the free Klein-Gordon equation is examined. The treatment in these papers is non-generic in the sense that the different branches of the dispersion relations of the systems considered are either polynomial or rational in the wave number (as in \cite{fokaspellonisys, pp1}) or are conveniently parameterized to avoid radicals, as in \cite{pp3, pp2, treharnefokas}. In addition, all five papers start from a Lax pair and construct a Riemann-Hilbert problem almost exclusively, with almost no hints as to how to generalize the more accessible UTM for scalar systems based on the use of Green's Theorem. Generically, the different branches of the dispersion relations for a first-order evolution system of dimension $N$ depend on radicals as they are the roots of an $N$-th order polynomial. Below we investigate the general case of an $N$-dimensional first-order linear system, using two examples. 

The first example is the (free: no potential) Klein-Gordon (KG) equation \cite{mandlshaw, pp3, pp1}

\beq\la{kg}
u_{tt}-u_{xx}+\alpha u=0,
\eeq

\no where $\alpha$ is a constant parameter and indices denote partial derivatives. For our purposes, the KG equation is rewritten as a two-dimensional first-order system:

 \begin{alpheqn}
\bea\la{fn}
q_{t}&=&p,\\
p_{t}&=&q_{xx}-\alpha q,
\eea
\end{alpheqn}

\no where $q=u$, $p=u_t$. Our approach to the KG equation differs from \cite{pp3, pp1} in that all our considerations are based on so-called local relations, to which Green's Theorem is applied. No parameterization of the dispersion relation branches is used. Rather we use the branches of the dispersion relation is their original form. Lastly, in \cite{pp3,pp2} only $\alpha>0$ is considered.

The second example is the linearized Fitzhugh-Nagumo (FN) system of partial differential equations \cite{fn}

\begin{alpheqn}
\bea\la{fn}
v_t&=&v_{xx}-v-w,\\
w_t&=&\beta v,
\eea
\end{alpheqn}

\no where $\beta$ is a constant parameter. In general, we consider systems of the form

\beq\la{gensys}
Q_t+\Lambda(-i \p_x)Q=0,
\eeq

\no where $Q$ is an $N$-dimensional vector and $\Lambda$ is a linear-operator matrix of size $N\times N$ and of order $n$.

\vs

{\bf Example 1.} For the KG equation, we have
\beq
\p_t \left(\ba{c}q\\p\ea\right)+\left(\ba{cc}0&-1\\-\p_x^2+\alpha&0\ea\right)\left(\ba{c}q\\p\ea\right)=0,
\eeq

\no so that $n=2$, and

\beq
\Lambda(-i\p_x)=\left(\ba{cc}0&-1\\(-i\p_x)^2+\alpha&0\ea\right).
\eeq

\vs

{\bf Example 2.} Similarly, for the FN system,
\beq
\p_t \left(\ba{c}v\\w\ea\right)+\left(\ba{cc}-\p_x^2+1&1\\-\beta&0\ea\right)\left(\ba{c}v\\w\ea\right)=0,
\eeq
allowing us to read off $\Lambda(-i\p_x)$ and $n=2$ as well.

\vs

Throughout, we contrast the general systems case with the general scalar case, where the equation can be written in the form

\beq\la{scalarpde}
u_t+\lambda(-i \p_x)u=0,
\eeq

\no for a scalar-valued function $u(x,t)$.

Although different nonlocal systems can be considered by allowing $\Lambda$ to depend rationally on its argument (see {\em e.g.} \cite{vd} for the scalar setting), in this paper we restrict to the case where $\Lambda$ (and $\lambda$) depends polynomially on its argument. As a consequence, the dispersion relation for the scalar case, easily found by equating $u=\exp(ikx-\omega t)$, is given by $\omega=\lambda(k)$, and $\omega$ depends polynomially on $k$. It should be remarked that in the above calculation of the dispersion relation, we have followed the convention in the literature on the UTM for the dispersion relation. Thus the dispersion relation as used here differs by a factor of $i$ from the standard use. For instance, in the UTM dispersive equations are characterized by a purely imaginary dispersion relation for $k\in \mathbb{R}$.

Similarly, for the systems case \rf{gensys}, we let

\beq
Q=\left(\ba{c}Q_1\\\vdots\\Q_N\ea\right)e^{ikx-\omega t},
\eeq

\no so that $\omega$ satisfies

\beq\la{dr}
\det(\Lambda(k)-\omega I)=0,
\eeq

\no where $I$ is the $N\times N$ identity matrix. Thus the different branches $\Omega_1$, \ldots, $\Omega_N$ of the dispersion relation are roots of an $N$-th order polynomial. Generically, they depend on radicals of order up to $N$ ($N$-th roots) whose arguments are polynomials in $k$. It follows that these dispersion branches are sheets of an $N$-valued function and they have branch point singularities in the complex $k$ plane. This is not always the case, as is illustrated in the appendix for the wave equation and in the examples treated in \cite{fokaspellonisys}. Throughout this paper, we assume that all branches $\Omega_j(k)$ are distinct, except at isolated values of $k\in \mathbb{C}$.

\vs

{\bf Example 1.} For the KG equation,
\beq
\Omega_{1,2}=\pm i \sqrt{\alpha+k^2},
\eeq
which has branch points at the square roots of $-\alpha$. We choose the branch cut that connects these branch points.

\vs

{\bf Example 2.} For the FN system,
\beq\la{fnbranches}
\Omega_{1,2}=\frac{1+k^2\pm\sqrt{(1+k^2)^2-4\beta}}{2},
\eeq
which has 4 branch points.

\vs

One of the main advantages of the UTM is its ability to characterize the number and the type of boundary conditions required to ensure wellposedness of a given initial-boundary value problem \cite{fokasbook}. We wish to see to what extent the same can be done for systems of PDEs. To this end, we do not specify boundary conditions for our examples at this point. Rather, we will see how the application of Fokas's UTM determines the information that should be provided on the boundary of our domain. We limit ourselves to problems posed on the half line $x>0$. It is anticipated that extending our results to problems on the finite interval $x\in (0,L)$ is comparable to extending the UTM for scalar problems on the half line to scalar problems on the finite interval. Of course, in all cases, we specify initial conditions $Q(x,0)=Q_0(x)$.

In the next sections, we go through the extension of Fokas's UTM, as applied to systems of linear, constant-coefficient evolution PDEs. Each section deals with a different step of the method, so as to present the method in an algorithmic way.

\section{The local relation}

The first step for the application of the UTM is to rewrite the system of equations in divergence form. We refer to this form as the {\bf local relation}. Following the appendix of \cite{fokaspellonisys}, we rewrite \rf{gensys} as

\begin{align}\la{localrelation}
&&\left(e^{-ikxI+\Lambda(k) t}Q\right)_t-\left(e^{-ikxI+\Lambda(k) t}X(x,t,k)Q\right)_x&=0\\\la{localrelationbranched}
&\Ra&\left(e^{-ikxI+\Omega(k) t}A(k)Q\right)_t-\left(e^{-ikxI+\Omega(k) t}A(k)X(x,t,k)Q\right)_x&=0,
\end{align}

\no where $\Omega(k)=\mbox{diag}(\Omega_1, \ldots, \Omega_N)$, the diagonal matrix with the different branches of the dispersion relation as diagonal elements. The matrix $A(k)$ diagonalizes the matrix $\Lambda(k)$:

\beq
\Lambda(k)=A^{-1}(k)\Omega(k)A(k).
\eeq

\no Lastly, the vector $X(x,t,k)$ is a differential matrix operator of degree at most $n-1$, polynomial in $k$, defined by

\beq
X(x,t,k)=i\left.\frac{\Lambda(k)-\Lambda(l)}{k-l}\right|_{l=-i\p_x}=\sum_{j=0}^{n-1} c_j(k)\p_x^j,
\eeq

\no and the last equality defines the matrix-valued polynomials $c_j(k)$. Equation \rf{localrelation} is verified by working out the product rules of both terms. For the scalar case, both equations are identical. For the case of systems, an important difference between \rf{localrelation} and \rf{localrelationbranched} is that \rf{localrelationbranched} contains $\Omega(k)$ and $A(k)$, which are typically branched in the complex $k$ plane. On the other hand, the left-hand side of \rf{localrelation} depends only on $\Lambda(k)$ and $X(x,t,k)$, which are not branched as functions of $k$. In practice, \rf{localrelationbranched} is more useful, as \rf{localrelation} requires the calculation of a matrix exponential with non-diagonal exponent. Even in specific examples where the matrix exponential is easily calculated, it is useful to have the local relation in terms of the branches of the dispersion relation.

\vs

{\bf Example 1.} For the KG equation, $\exp(\Lambda(k)t)$ is easily computed directly:

\beq
e^{\Lambda(k) t}=I \cos(\sqrt{\alpha+k^2}t)+\Lambda(k)\,\mbox{sinc}(\sqrt{\alpha+k^2}t),
\eeq

\no which is not branched since all functions above with square root arguments are even. With

\beq
A=\twomatrix{-\Omega_2(k)}{-1}{\Omega_1(k)}{1}, ~~X\twovector{q}{p}=\twovector{0}{ikq+q_x},
\eeq

\no the local relations \rf{localrelationbranched} are

\beq\la{KGlocrel}
\left(e^{-ikx+\Omega_j(k)t}(-\Omega_j(k)q+p)\right)_t-\left(e^{-ikx+\Omega_j(k)t}(ikq+q_x)\right)_x=0,~~~j=1,2.
\eeq

\vs

{\bf Example 2.} For the FN equation,
\beq
A=\twomatrix{-\beta}{-\Omega_2(k)}{\beta}{\Omega_1(k)}, ~~X\twovector{v}{w}=\twovector{ikv+v_x}{0},
\eeq

\no and the local relations \rf{localrelationbranched} are

\beq\la{FNlocrel}
\left(e^{-ikx+\Omega_j(k)t}(\Omega_j(k)v+w)\right)_t-\left(e^{-ikx+\Omega_j(k)t}(ik\Omega_j(k)v+\Omega_j(k)v_x)\right)_x=0,
~~~j=1,2,
\eeq
after multiplication by $\Omega_j(k)$, and using that $\Omega_1(k)\Omega_2(k)=\beta$.

\vs

{\bf Remarks.}
\begin{itemize}

\item For the scalar case, there is no difference between the two forms \rf{localrelation} and \rf{localrelationbranched}. The local relation is easily obtained by multiplying the PDE \rf{scalarpde} by $\exp(-ikx+\omega t)$ and using integration by parts to get to the divergence form. The dispersion relation is found during this process \cite{dtv} as well. The same process can be used for systems of PDEs, but additional linear algebra is required to get to the divergence form, as only specific linear combinations of the equations allow for this form.

\item In practice, it may not always be possible to write the set of local relations as compactly as above, using index notation: the first local relation may depend on $\Omega_2(k)$ and so on. In fact, this is the case in the second example, but the use of the dispersion relation for the FN system allows the further simplification.

\item In \cite{treharnefokas} and to a lesser extent in \cite{fokaspellonisys}, the Lax pair formalism for the UTM is used. For constant-coefficient systems linear PDEs this is not necessary, and we stay within the framework of using the local relation and Green's Theorem (in the next step).

\end{itemize}

\section{The global relation}

\begin{figure}
\centering
\includegraphics[width=3in]{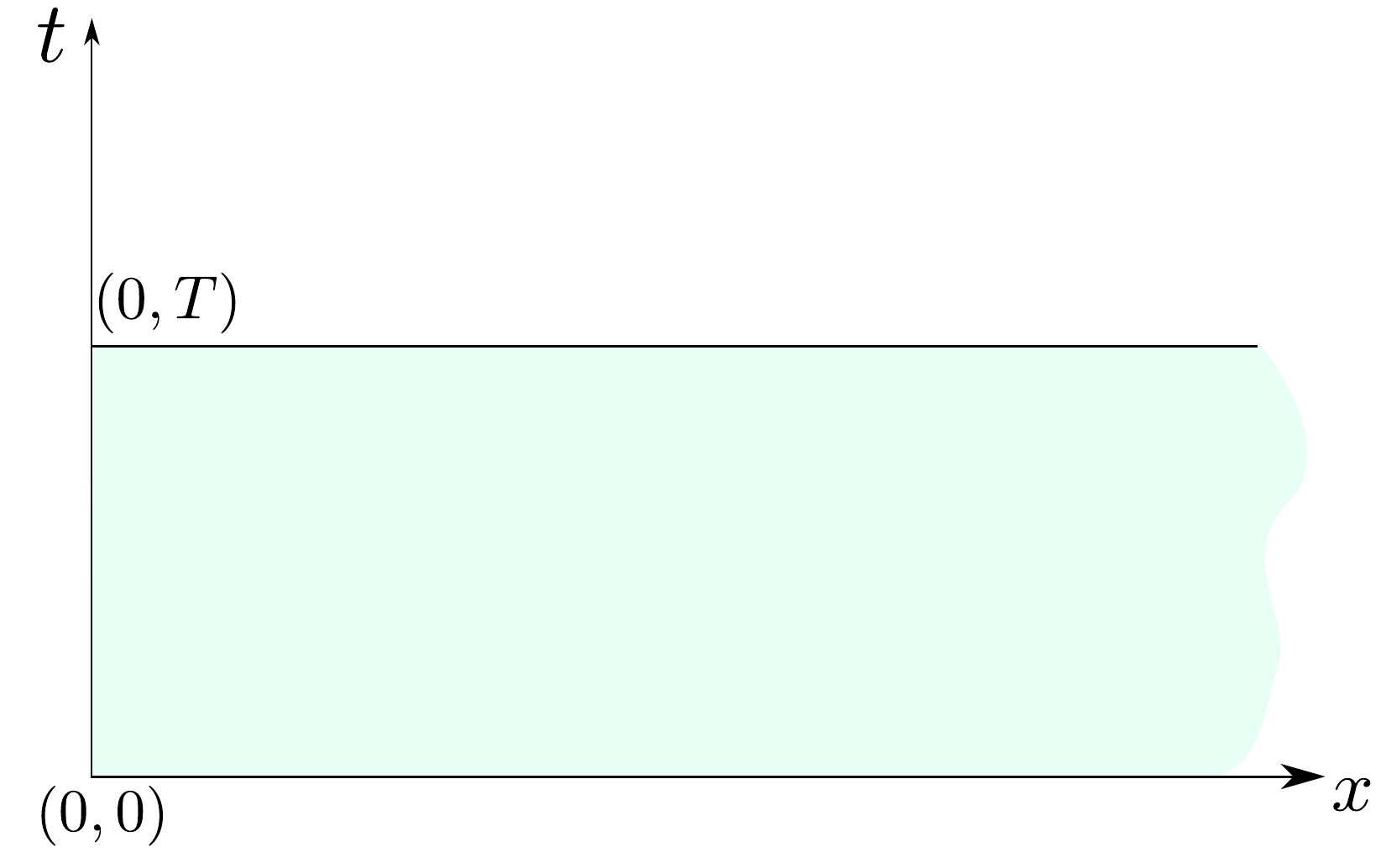}
\caption{The region of integration in the $(x,t)$ plane for boundary-value problems posed on the positive half-line.}\label{fig:half-line}
\end{figure}

For the next step, we integrate each local relation over an infinite strip in the $(x,t)$ plane, cornered at the origin, see Fig.~\ref{fig:half-line}. Using Green's Theorem, we obtain

\beq\la{globalrelation}
\hat Q_0(k)-e^{\Lambda(k)t}\hat Q(k,t)-G(k,t)=0,
\eeq

\no using \rf{localrelation}, and

\beq\la{globalrelationbranched}
A(k)\hat Q_0(k)-e^{\Omega(k)t}A(k)\hat Q(k,t)-\tilde G(k,t)=0,
\eeq

\no using \rf{localrelationbranched}. Here

\alpheqn
\begin{align}\la{defs}
\hat Q_0(k)=\int_{0}^\infty e^{-ikx} Q_0(x)dx, && \hat Q(k,t)=\int_{0}^\infty e^{-ikx} Q_(x,t)dx,\\\la{defsb}
G(k,t)=\int_0^t e^{\Lambda s}X(0,s,k)Q(0,s)ds, && \tilde G(k,t)=\int_0^t e^{\Omega s}A(k)X(0,s,k)Q(0,s)ds,
\end{align}
\resetalpheqn

Note that we have replaced the upper limit $T$ (see Fig.~\ref{fig:half-line}) of the temporal variable by $t$. The relations \rf{globalrelation} and \rf{globalrelationbranched} are referred to as the {\bf Global Relations}. The same comments can be made as for the local relations: \rf{globalrelation} is written entirely in terms of quantities without branch points, while \rf{globalrelationbranched} contains branched quantities throughout, but is easier to write down in practice. For specific examples, it is more convenient to start directly from the derived local relations.

{\bf Example 1.} For the KG equation,

\beq\la{grkg}
-\Omega_j \hat q_0+\hat p_0-e^{\Omega_j t}(-\Omega_j \hat q+\hat p)-ik g_0(\Omega_j,t)-g_1(\Omega_j,t)=0, ~~j=1,2,
\eeq

\no where
\beq
\hat f_0(k)=\int_0^\infty e^{-ikx} f_0(x)dx, ~~\hat f(k,t)=\int_0^\infty e^{-ikx} f(x,t)dx,~~
g_\kappa(\Omega_j,t)=\int_0^t e^{\Omega_j s} q_{\kappa x}(0,s)ds,
\eeq

\no for $f=q,~p$, $\kappa=0,1$ and the index $\kappa x$ denotes $\kappa$ derivatives with respect to $x$.

\vs

{\bf Example 2.} For the FN system,

\beq\la{grfn}
\Omega_j \hat v_0+\hat w_0-e^{\Omega_j t}(\Omega_j \hat v+\hat w)-ik \Omega_j g_0(\Omega_j,t)-\Omega_j g_1(\Omega_j,t)=0, ~~j=1,2,
\eeq

\no where the hat denotes the half-line Fourier transform, as above, and
\beq
g_\kappa(\Omega_j,t)=\int_0^t e^{\Omega_j s} v_{\kappa x}(0,s)ds, ~~\kappa=0, 1.
\eeq

\vs

As in the case of constant-coefficient scalar evolution equations posed on the half line, the global relation is valid in the lower half of the complex plane, $\Im\, k\leq 0$. Indeed, the half-line Fourier transforms in \rf{defs} are defined only if the exponential does not grow, assuming integrable or square-integrable initial conditions and solutions. The other quantities are defined by proper integrals and do not impose any further restriction.

\section{A ``solution'' formula}

The global relations contain the Fourier transform of the quantities we wish to solve for. At this point, we solve for these transforms and invert using a regular Fourier inversion, interpreting the half-line transforms as whole-line transforms of functions that are zero for $x<0$. Thus

\beq\la{sol}
Q(x,t)=\frac{1}{2\pi} \int_{-\infty}^\infty e^{ikxI-\Lambda(k)t}\hat Q_0(k)dk-\frac{1}{2\pi}\int_{-\infty}^\infty e^{ikxI-\Lambda(k)t}G(k,t)dk,
\eeq

\no using \rf{globalrelation}, and, clearly equivalent,

\beq\la{solbranched}
Q(x,t)=\frac{1}{2\pi} \int_{-\infty}^\infty A^{-1}(k)e^{ikxI-\Omega(k)t}A(k)\hat Q_0(k)dk-\frac{1}{2\pi}\int_{-\infty}^\infty A^{-1}(k) e^{ikxI-\Omega(k)t} \tilde G(k,t)dk.
\eeq

These equations do not provide expressions for the solutions of \rf{gensys}. Although the first term contains only known quantities, the second term depends on $m_j-1$ derivatives of $Q_j$ evaluated at the boundary $x=0$. Here $Q_j$ is the $j$-th component of $Q$, $j=1, \ldots, N$ and $m_j$ is the maximal degree as a function of $k$ occurring in the $j$-th column of $\Lambda(k)$. If $m_j=0$, the function $Q_j$ nor any of its $x$ derivatives appear in the second term and no boundary data involving $Q_j$ should be prescribed. We know from experience with scalar problems \cite{dtv, fokasbook} that it is unlikely that all these unknown boundary functions need to be specified. In fact, specifying them freely is likely to result in inconsistencies: {\em e.g.}, it is not possible to prescribe the Dirichlet and Neumann data for the heat equation on the half line independently.

\vs

{\bf Example 1.} For the KG equation,
\alpheqn
\begin{align}\nonumber
q(x,t)=&\frac{1}{2\pi}\int_{-\infty}^\infty \frac{e^{ikx}}{\Omega_1-\Omega_2}\left(
e^{-\Omega_1 t}(\Omega_1 \hat q_0-\hat p_0)-e^{-\Omega_2 t}(\Omega_2 \hat q_0-\hat p_0)
\right)dk\\\la{kgsol}
&-\frac{1}{2\pi}\int_{-\infty}^\infty \frac{e^{ikx}}{\Omega_1-\Omega_2}\left(
e^{-\Omega_1 t}(-ik g_0^{(1)}-g_1^{(1)})-
e^{-\Omega_2 t}(-ik g_0^{(2)}-g_1^{(2)}
\right)dk,\\\nonumber
p(x,t)=&\frac{1}{2\pi}\int_{-\infty}^\infty \frac{e^{ikx}}{\Omega_1-\Omega_2}\left(
\Omega_2 e^{-\Omega_1 t}(\Omega_1 \hat q_0-\hat p_0)-\Omega_1 e^{-\Omega_2 t}(\Omega_2 \hat q_0-\hat p_0)
\right)dk\\
&-\frac{1}{2\pi}\int_{-\infty}^\infty \frac{e^{ikx}}{\Omega_1-\Omega_2}\left(
\Omega_2 e^{-\Omega_1 t}(-ik g_0^{(1)}-g_1^{(1)})-\Omega_1
e^{-\Omega_2 t}(-ik g_0^{(2)}-g_1^{(2)}
\right)dk,
\end{align}
\resetalpheqn

\no where $g_\kappa^{(j)}=g_\kappa(\Omega_j,t)$, $\kappa=0,1$, $j=1,2$.

\vs

{\bf Example 2.} For the FN system,
\alpheqn
\begin{align}\nonumber
v(x,t)=&\frac{1}{2\pi\beta}\int_{-\infty}^\infty \frac{e^{ikx}}{\Omega_2-\Omega_1}\left(
\Omega_2 e^{-\Omega_2 t}(\beta \hat v_0+\Omega_1 \hat w_0)-\Omega_1 e^{-\Omega_1 t}(\beta \hat v_0+\Omega_2 \hat w_0)
\right)dk\\\la{fnsol}
&-\frac{1}{2\pi}\int_{-\infty}^\infty \frac{e^{ikx}}{\Omega_2-\Omega_1}\left(
\Omega_2 e^{-\Omega_2 t}(ikg_0^{(2)}+g_1^{(2)})-\Omega_1 e^{-\Omega_1 t}(ikg_0^{(1)}+g_1^{(1)})
\right)dk,\\\nonumber
w(x,t)=&\frac{1}{2\pi}\int_{-\infty}^\infty \frac{e^{ikx}}{\Omega_2-\Omega_1}\left(
e^{-\Omega_1 t}(\beta \hat v_0+\Omega_2 \hat w_0)-e^{-\Omega_2 t}(\beta \hat v_0+\Omega_1 \hat w_0)
\right)dk\\
&-\frac{\beta}{2\pi}\int_{-\infty}^\infty \frac{e^{ikx}}{\Omega_2-\Omega_1}\left(
e^{-\Omega_1 t}(ikg_0^{(1)}+g_1^{(1)})- e^{-\Omega_2 t}(ikg_0^{(2)}+g_1^{(2)})
\right)dk,
\end{align}
\resetalpheqn
\no using the same notation as above.

\vs

In the following sections, we aim to turn \rf{sol} and \rf{solbranched} into genuine solution formulas, depending on the correct number of boundary functions.

\section{Deformation of the integration path}

Following the approach to solve constant-coefficient scalar equations, we wish to deform the integration path in the expressions \rf{sol} and \rf{solbranched} as far away from the real line as possible as $k\ra \infty$ \cite{dtv, fokasbook}. In some cases all or part of the real line will remain in place, and other curves in the upper-half plane resulting in zero contributions might be added to the integration path, as in the scalar case. As in that case, we have no intent to deform the path of integration of the first term: the integrand of this term is known explicitly and it requires no further manipulation.

Thus we focus on the second term. Specifically, we wish to elucidate the role of the branch points in the systems case. Since \rf{sol} is written in terms of quantities that are not branched, it is clear that the branch points of $A(k)$ and $\Omega(k)$ that are apparent in \rf{solbranched} are all removable. This is also immediately clear from the expressions (\ref{kgsol}-b) and (\ref{fnsol}-b). Encircling any of the branch points results in an interchange of the two sheets, effectively switching the indices on $\Omega_1$ and $\Omega_2$. This permutation leaves the integrands invariant, confirming the removability of their branch points. As stated above, this is true in general, due to the non-branched nature of the integrands in \rf{sol}. For any given example, the ``solution formula'' is conveniently written using \rf{solbranched}, resulting in integrands that are symmetric functions of the different branches of the dispersion relation. Encircling a branch point amounts to a permutation of the indices of these branches, producing no change due to the symmetry of the integrands. It follows that the process of deforming the path of integration proceeds very much as in the scalar case, with minor modifications due to the vector structure as discussed below.

We examine the $k$ dependence of the integrand in the second term of \rf{sol}. Since $x$ and $t$ are independent variables, both defined on their respective half lines, different restrictions on $k\in \mathbb{C}$ are imposed by the need to control the $x$ and $t$ dependence of the integrand separately. Since the $x$ dependence is confined to the exponential, it follows that necessarily $k\in \mathbb{C}^+=\{\mathbb{C}: \Im\,k\geq 0\}$ for the integral to be defined. Similarly, using \rf{defsb} and following the reasoning from the scalar case \cite{dtv, fokasbook}, we need $\Re\, \Omega(k)\geq 0$, {\em i.e.}, the matrix $\Omega(k)$ is positive definite. Thus, we define the inaccessible region $D^+$ in the upper-half plane:

\beq\la{badregion}
D^+=\bigcup_{j=1}^N\{k\in \mathbb{C}: \Im\,k>0,~\Re\,\Omega_j(k)<0\}.
\eeq

\no Further, for $k\in \mathbb{C}^+/D^+$, the integrand decays exponentially as $k\ra \infty$. Thus the contribution to the second term of any path that tends to infinity in $\mathbb{C}^+/D^+$ vanishes, by Jordan's Lemma \cite{af}. It follows from Cauchy's Theorem \cite{af} that we may write

\beq\la{soldef}
Q(x,t)=\frac{1}{2\pi} \int_{-\infty}^\infty e^{ikxI-\Lambda(k)t}\hat Q_0(k)dk-\frac{1}{2\pi}\int_{\p D^+} e^{ikxI-\Lambda(k)t}G(k,t)dk,
\eeq

\no using \rf{globalrelation}, and, clearly equivalent,

\beq\la{solbrancheddef}
Q(x,t)=\frac{1}{2\pi} \int_{-\infty}^\infty A^{-1}(k)e^{ikxI-\Omega(k)t}A(k)\hat Q_0(k)dk-\frac{1}{2\pi}\int_{\p D^+} A^{-1}(k) e^{ikxI-\Omega(k)t} \tilde G(k,t)dk.
\eeq

\vs

\sloppypar The integrand of the second term is a linear combination of different exponentials $\exp{(ikx-\Omega_j t)}$. One may consider splitting the integral into different parts, each depending on one of these exponentials only. If this is done, the different integrands are not symmetric functions of the branched quantities, and the branch points are not removable. As a consequence, much greater care is required for the deformation of the integration path away from the real line. To facilitate this, it is convenient to deform the integration path above all finite and possibly infinite branch points of $\Omega(k)$, see Fig.~\ref{fig:deform1}. The deformation above the finite (removable) branch points is trivially allowed by Cauchy's Theorem. Deforming around branch points at infinity should be considered on a case-by-case basis, as convergence issues come into play. We will not consider this case further.

After the initial finite deformation, further deformations into the region not containing the branch points can be made without these points playing any role at all. Once the branch points are below the path of integration, it may be possible and desirable that different integration paths $\p D^+_j$ are used for different components $Q_j$, if these components do not depend on all $\Omega_j$. 


\begin{figure}[tb]
\begin{center}
\def\svgwidth{5in}
    \vspace*{0in}
    \hspace*{0.0in}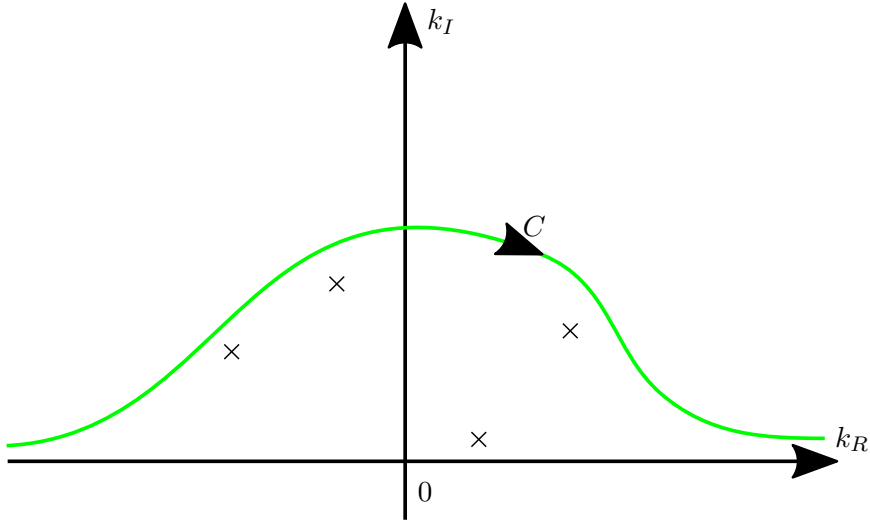
    \vspace*{0in}
    \caption{\la{fig:deform1} Deforming away from the real line, above all finite branch points in $\mathbb{C}^+$, with $k=k_R+i k_I$.}
\end{center}
\end{figure}

\vs

{\bf Example 1.} For the KG equation $q$ and $p$ depend on both $\Omega_1$ and $\Omega_2$. In general, $D^+$ is conveniently found by first determining its boundary, followed by examining the different regions it defines. For the KG equation, \rf{dr} reads

\beq
\omega^2+\alpha+k^2=0.
\eeq

\no Splitting this expression in its real and imaginary parts and imposing that Re\,$\omega=0$ determines the boundary of $D^+$. This boundary is shown in red in Figure~\ref{fig:kgdeform}. This boundary does not separate the upper-half plane in different regions. Since $\Omega_1=-\Omega_2$, there is no subset of the upper-half plane where the real part of both $\Omega_1$ and $\Omega_2$ is positive. Thus $D^+$ consists of the whole upper-half plane and no deformation from the integration along the real line to the upper-half plane in (\ref{kgsol}-b) is done.

\vs

\no {\bf Remark.} Alternatively, one may deform the integration to go above the branchpoint at $i\sqrt{\alpha}$ ($\alpha>0$) or those at $\pm\sqrt{-\alpha}$ ($\alpha<0$). After doing so, the integral on the second term may be distributed, resulting in one integral where the whole upper-half plane is accessible. By Cauchy's Theorem, the contribution from this integral vanishes. In other words, the net contribution from this branch of the dispersion relation is hiding in the other integral as the contribution from the branch point(s) as the other integral is deformed back to the real line.

\begin{figure}[tb]
\begin{center}
\begin{tabular}{cc}
\def\svgwidth{2.8in}
    \vspace*{0in}
    \hspace*{0.0in}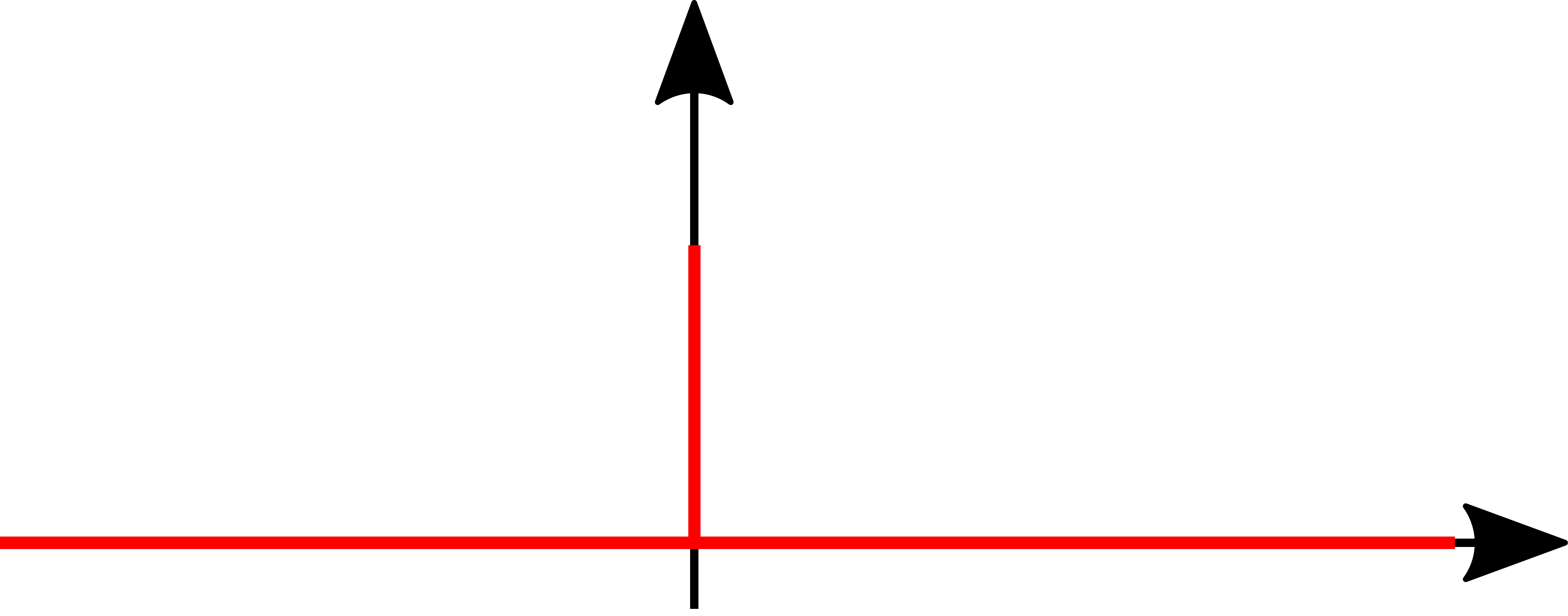
    \vspace*{0in}&
    \def\svgwidth{2.8in}
    \vspace*{0in}
    \hspace*{0.0in}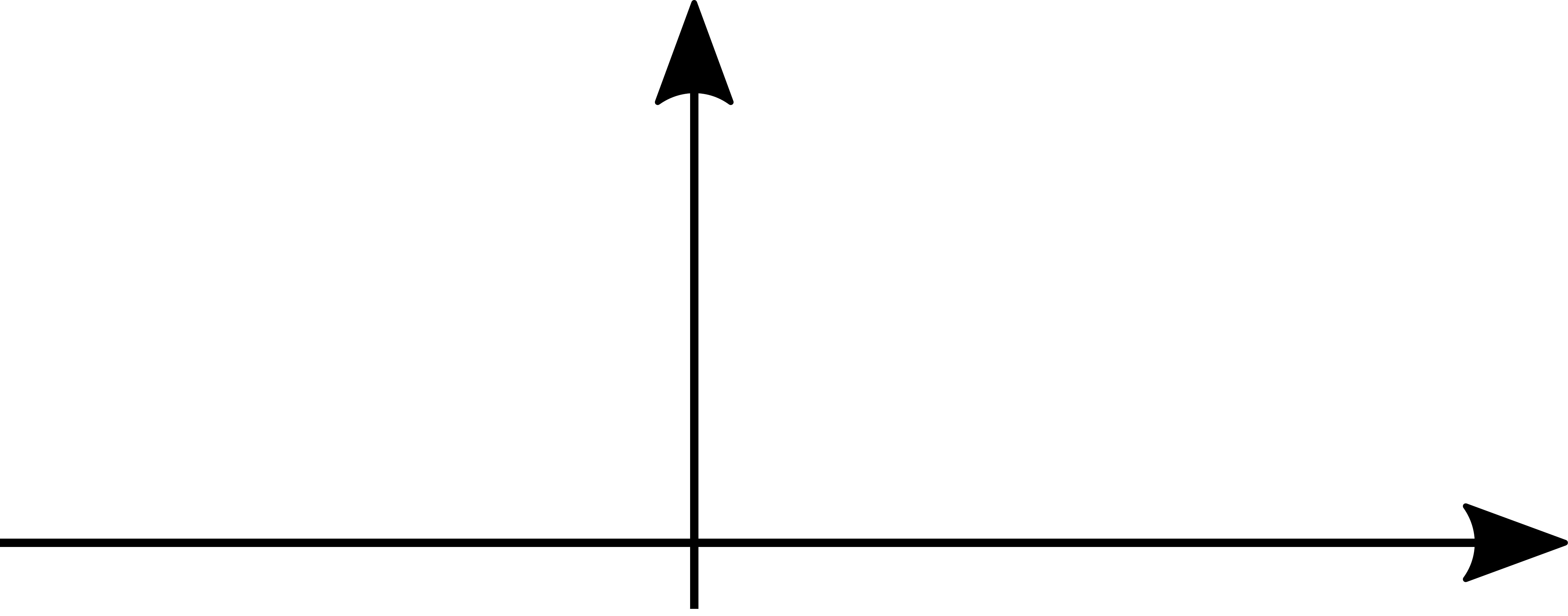
    \vspace*{0in}
\end{tabular}
    \caption{\la{fig:kgdeform} The region $D+$, with boundary shown in red, covering the whole upper-half plane, for $\alpha>0$ (left) and $\alpha<0$ (right).}
\end{center}
\end{figure}

\vs

{\bf Example 2.} For the FN equation the situation is more straightforward.
\beq
D^+=\{k=k_R+i k_I\in \mathbb{C}: k_I>0~\mbox{and}~k_I^2>1+k_R^2\},
\eeq

\no which is illustrated in Fig.~\ref{fig:fnregion}. Depending on the value of $\beta$, some of the different branch points of $\Omega_1$ and $\Omega_2$ may lie between the real line and $\p D^+$. As long as the integrals are not distributed over their integrands, this is not a concern. Thus

\alpheqn
\begin{align}\nonumber
v(x,t)=&\frac{1}{2\pi\beta}\int_{-\infty}^\infty \frac{e^{ikx}}{\Omega_2-\Omega_1}\left(
\Omega_2 e^{-\Omega_2 t}(\beta \hat v_0+\Omega_1 \hat w_0)-\Omega_1 e^{-\Omega_1 t}(\beta \hat v_0+\Omega_2 \hat w_0)
\right)dk\\\la{fnsoldef}
&-\frac{1}{2\pi}\int_{\p D^+} \frac{e^{ikx}}{\Omega_2-\Omega_1}\left(
\Omega_2 e^{-\Omega_2 t}(ikg_0^{(2)}+g_1^{(2)})-\Omega_1 e^{-\Omega_1 t}(ikg_0^{(1)}+g_1^{(1)})
\right)dk,\\\nonumber
w(x,t)=&\frac{1}{2\pi}\int_{-\infty}^\infty \frac{e^{ikx}}{\Omega_2-\Omega_1}\left(
e^{-\Omega_1 t}(\beta \hat v_0+\Omega_2 \hat w_0)-e^{-\Omega_2 t}(\beta \hat v_0+\Omega_1 \hat w_0)
\right)dk\\
&-\frac{\beta}{2\pi}\int_{\p D^+} \frac{e^{ikx}}{\Omega_2-\Omega_1}\left(
e^{-\Omega_1 t}(ikg_0^{(1)}+g_1^{(1)})- e^{-\Omega_2 t}(ikg_0^{(2)}+g_1^{(2)})
\right)dk.
\end{align}
\resetalpheqn

\begin{figure}[tb]
\begin{center}
\def\svgwidth{5in}
    \vspace*{0in}
    \hspace*{0.0in}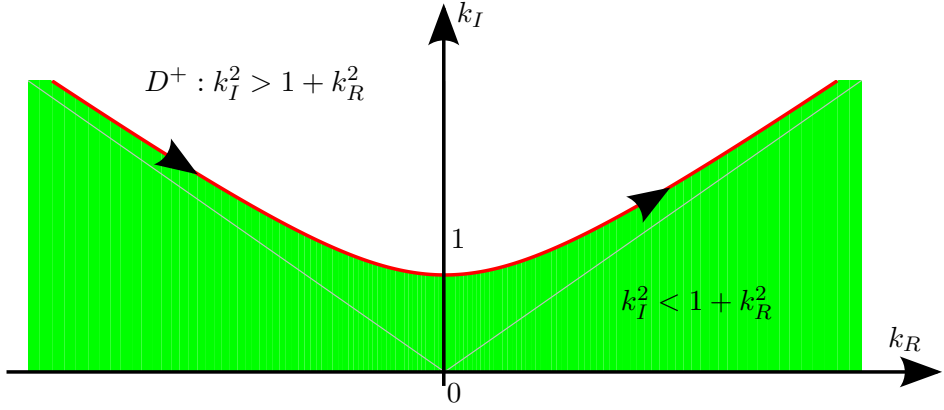
    \vspace*{0in}
    \caption{\la{fig:fnregion} The region $D^+$ for the FN equation and the deformed path of integration along its boundary in red.}
\end{center}
\end{figure}

\section{Symmetries of the dispersion relation}

In order to eliminate the superfluous dependence on unknown boundary conditions, we need the symmetry group of the dispersion relation. This is the collection of transformations $k\ra \nu(k)$ that leave the dispersion relation invariant. Applying these transformations allows us to transform one branch of the dispersion relation into another one. The symmetries are found by solving

\beq\la{symmetries}
\det (\Lambda(\nu)-\Omega_j(k)I)=0,~~~j=1, \ldots, N.
\eeq

\no which determines $\nu$ as a function of $k$, since $\Omega_j$ depends on $k$ through \rf{dr}. It is clear that the identity transformation $\nu_1=k$ is a solution of \rf{symmetries}. Below we show that for both the KG and FN system, the only nontrivial symmetry is $\nu_2: k\ra -k$. The appendix contains an example of a different symmetry.

\vs

{\bf Example 1.} For the KG equation, \rf{symmetries} is
\begin{align}\nonumber
&&\nu^2+\alpha+\Omega_j^2&=0\\\nonumber
&\Ra& \nu^2-k^2&=0\\
&\Ra& \nu_1=k, ~~&\nu_2=-k,
\end{align}
for $j=1, 2$.

\vs

{\bf Example 2.} For the FN equation, for $j=1,2$,
\begin{align}\nonumber
&&\Omega_j^2-(1+\nu^2)\Omega_j+\beta&=0\\\nonumber
&\Ra& (1+k^2)\Omega_j-\beta-(1+\nu^2)\Omega_j+\beta&=0\\\nonumber
&\Ra& (k^2-\nu^2)\Omega_j&=0\\
&\Ra& \nu_1=k, ~~&\nu_2=-k,
\end{align}
where we have used that $\Omega_j^2=(1+k^2)\Omega_j-\beta$.

\section{The elimination of unknown boundary functions}

The deformed solution formula \rf{solbrancheddef} depends on the quantities
\beq
g_{j,l}^{(m)}=g_{j,l}(\Omega_m,t)=\int_0^t e^{-\Omega_m s}Q_{j,lx}(0,s)ds,
\eeq
where $Q_{j,lx}$ denotes $l$ derivatives with respect to $x$ of $Q_j(x,t)$. Through these quantities, the solution formula exhibits dependence on a number of boundary functions $Q_{j,lx}$. In order to eliminate such unnecessary boundary functions from the deformed solution formulae \rf{solbrancheddef}, new global relations \rf{globalrelationbranched} are obtained by using the symmetries from the previous section.

Since the matrix $\exp(\Omega s)$ is diagonal, the $m$-th component of \rf{globalrelationbranched} contains only $g_{j,k}^{(m)}$, for varying $j$ and $k$. In other words, $\Omega_m$ is the only branch of the dispersion relation appearing in the $m$-th component. If $\Omega_m(\nu_l(k))=\Omega_n(k)\in \{\Omega_j,~j=1, \ldots, N\}$, the new global relation below is considered:
\beq\la{ngr}
\sum_{j=1^N}A_{mj}(\nu_l(k))\hat Q_{0,j}(\nu_l(k))-e^{\Omega_n(k)t}\sum_{j=1^N}A_{mj}(\nu_l(k))\hat Q_{j}(\nu_l(k))-\tilde G_m(\nu_l(k))=0,
\eeq
with
\beq
\tilde G_m(\nu_l(k))=\int_0^t e^{\Omega_n(k)s}\left(A(\nu_l(k))X(0,s,\nu_l(k))Q(0,s)\right)_m ds,
\eeq
where the subindex $m$ indicates the $m$-th component is taken. It follows that these new global relations depend on the same quantities $g_{j,l}^{(m)}$ as does the deformed solution formula \rf{solbrancheddef}. The new global relation is valid for all $k$ for which $\nu_l(k)\in \{k\in \mathbb{C}: \Im k\leq 0\}$. In other words, $k$ is in the closed lower-half plane transformed under $\nu_l^{-1}$. This process results in a number of new global relations. If these global relations are valid in regions of $\mathbb{C}$ where all or a part $\partial D^+$ lies, they can be used to eliminate unwanted boundary functions by solving the new global relations for these unwanted boundary conditions, as in the scalar case \cite{dtv, fokasbook}. The above procedure determines the exact number of boundary conditions that needs to be prescribed in order for the boundary-value problem to be well posed: this number is $\sum_{j=1}^N (m_j+1)$ minus the number of boundary functions that can be eliminated.

\vs

{\bf Example 1.} For the KG equation, there is only the nontrivial symmetry $\nu_2(k)=-k$. Both $\Omega_1(k)$ and $\Omega_2(k)$ are invariant under $k\ra -k$, thus we obtain two additional global relations, both valid in the closed upper-half plane. Importantly, both are valid on the real line, thus both can be used in \rf{kgsol}. The two new global relations are

\beq\la{grkgn}
-\Omega_j \hat q_0(-k)+\hat p_0(-k)-e^{\Omega_j t}(-\Omega_j \hat q(-k)+\hat p(-k))+ik g_0^{(j)}=g_1^{(j)}, ~~j=1,2.
\eeq

\no The solution formula \rf{kgsol} depends on $g_0^{(1)}$, $g_0^{(2)}$, $g_1^{(1)}$, and $g_1^{(2)}$. Using the two new global relations above, we expect to eliminate two of these. Assuming that Dirichlet boundary conditions are specified, we wish to eliminate the dependence in \rf{kgsol} on $g_1^{(1)}$ and $g_1^{(2)}$, which encode the Neumann data. Note that no boundary-value dependence on $p(x,t)$ arises. Using the expressions above, \rf{kgsol} becomes

\alpheqn
\begin{align}\nonumber
q(x,t)=&\frac{1}{2\pi}\int_{-\infty}^\infty
e^{ikx}(\hat q_0(k)-\hat q_0(-k))\frac{\Omega_1 e^{-\Omega_1 t}-\Omega_2 e^{-\Omega_2 t}}{\Omega_1-\Omega_2}dk+\\\nonumber
&-\frac{1}{2\pi}\int_{-\infty}^\infty
e^{ikx}(\hat p_0(k)-\hat p_0(-k))\frac{e^{-\Omega_1 t}- e^{-\Omega_2 t}}{\Omega_1-\Omega_2}dk+\\\la{kgfinal}
&+\frac{i}{\pi}\int_{-\infty}^\infty k e^{ikx}\frac{g_0^{(1)}e^{-\Omega_1 t}-g_0^{(2)}e^{-\Omega_2 t}}{\Omega_1-\Omega_2}dk+R_1,\\ \nonumber
p(x,t)=&\frac{1}{2\pi}\int_{-\infty}^\infty
e^{ikx}(\hat q_0(k)-\hat q_0(-k))\frac{\Omega_1\Omega_2 (e^{-\Omega_1 t}- e^{-\Omega_2 t})}{\Omega_1-\Omega_2}dk+\\\nonumber
&-\frac{1}{2\pi}\int_{-\infty}^\infty
e^{ikx}(\hat p_0(k)-\hat p_0(-k))\frac{\Omega_2 e^{-\Omega_1 t}- \Omega_1 e^{-\Omega_2 t}}{\Omega_1-\Omega_2}dk+\\
&+\frac{i}{\pi}\int_{-\infty}^\infty k e^{ikx}\frac{g_0^{(1)}\Omega_2 e^{-\Omega_1 t}-\Omega_1 g_0^{(2)}e^{-\Omega_2 t}}{\Omega_1-\Omega_2}dk+R_2.
\end{align}
\resetalpheqn

\no Here
\beq\la{r1r2}
R_1=\frac{1}{2\pi}\int_{-\infty}^\infty e^{ikx}\hat q(-k,t)dk, ~~R_2=-\frac{1}{\pi}\int_{-\infty}^\infty e^{ikx}\hat p(-k,t)dk,
\eeq

\no the contributions in the right-hand side of \rf{kgfinal} that exhibit dependence on the left-hand side. We have not used that $\Omega_1=-\Omega_2$ or that $\Omega_1^2=\Omega_2^2=-(\alpha+k^2)$, so as to exhibit the symmetry of the solution. It is a straightforward check from the above formulae that $q_t=p$.

\vs

{\bf Example 2.} As above, the FN equation only has the nontrivial symmetry $\nu_2(k)=-k$. Both $\Omega_1(k)$ and $\Omega_2(k)$ are invariant under $k\ra -k$, thus two additional global relations are obtained, both valid in the closed upper-half plane, which is where $\p D^+$ is. Thus both can be used in \rf{fnsoldef}. The two new global relations, solved for $g_0^{(j)}$ are

\beq\la{grkgn}
- \hat v_0(-k)-\frac{\hat w_0(-k)}{\Omega_j}+e^{\Omega_j t}\left(\hat v(-k)+\frac{\hat w(-k)}{\Omega_j}\right)+g_1^{(j)}=ik g_0^{(j)}, ~~j=1,2.
\eeq

\no The solution formula \rf{fnsol} depends on $g_0^{(1)}$, $g_0^{(2)}$, $g_1^{(1)}$, and $g_1^{(2)}$. Using the two new global relations above, we expect to eliminate two of these. Assuming that Neumann boundary conditions are specified, we wish to eliminate the dependence in \rf{kgsol} on $g_0^{(1)}$ and $g_0^{(2)}$, which encode the Dirichlet data. No boundary-value dependence on $w(x,t)$ arises. Using the expressions above, \rf{fnsoldef} becomes

\alpheqn
\begin{align}\nonumber
v(x,t)=&\frac{1}{2\pi}\int_{-\infty}^\infty e^{ikx}(\hat v_0(k)+\hat v_0(-k))\frac{\Omega_2 e^{-\Omega_2 t}-\Omega_1 e^{-\Omega_1 t}}{\Omega_2-\Omega_1}dk+ \\\nonumber
&+\frac{1}{2\pi}\int_{-\infty}^\infty e^ikx(\hat w_0(k)+\hat w_0(-k))\frac{e^{-\Omega_2 t}-e^{-\Omega_1 t}}{\Omega_2-\Omega_1}dk+ \\\la{fnfinal}
&-\frac{1}{\pi}\int_{\p D^+} e^{ikx}\frac{\Omega_2 e^{-\Omega_2 t} g_1^{(2)}-\Omega_1 e^{-\Omega_1 t} g_1^{(1)}}{\Omega_2-\Omega_1}dk+S_1,\\\nonumber
w(x,t)=& \frac{\beta}{2\pi}\int_{-\infty}^\infty e^{ikx} (\hat v_0(k)+\hat v_0(-k))\frac{e^{-\Omega_1 t}-e^{-\Omega_2 t}}{\Omega_2-\Omega_1}dk+\\\nonumber
&+\frac{1}{2\pi}\int_{-\infty}^\infty e^{ikx} (\hat w_0(k)+\hat w_0(-k))\frac{\Omega_2 e^{-\Omega_1 t}-\Omega_1 e^{-\Omega_2 t}}{\Omega_2-\Omega_1}dk+\\
&-\frac{\beta}{\pi}\int_{\p D^+} e^{ikx}\frac{g_1^{(1)}e^{-\Omega_1 t}-g_1^{(2)}e^{-\Omega_2 t}}{\Omega_1-\Omega_2}dk-S_2,
\end{align}
\resetalpheqn

\no where
\beq\la{s1s2}
S_1=-\frac{1}{2\pi}\int_{\p D^+} e^{ikx} \hat v(-k,t)dk, ~~S_2=-\frac{1}{2\pi}\int_{\p D^+} e^{ikx} \hat w(-k,t)dk.
\eeq

\no In \rf{fnfinal}, the path of integration along $\p D^+$ for the terms involving initial conditions has been deformed back to the real line, to allow their combination with the initial-condition terms already present in \rf{fnsoldef}.

\vs

In the above examples $R_1$, $R_2$, $S_1$ and $S_2$ represent the right-hand terms depending on the solution we wish to obtain. In the next section, we show these terms are zero.

\section{A solution formula}

The new global relations \rf{ngr} depend on the solution through $\hat Q_j(\nu_l(k))$. Thus solving these relations for boundary functions and substitution in the solution formula \rf{solbrancheddef} introduces $\hat Q_j(\nu_l(k))$ in the right-hand side of \rf{solbrancheddef}. It appears no effective solution formula has been obtained yet. At this point, we single out the right-hand terms depending on the solution, and we wish to show that their contribution is zero, as in the examples below. In general, in these culprit terms the time-dependent part of the exponential cancels, allowing the deformation of the contour into the previously inaccessible region. Further, for all examples we have examined, the right-hand side terms involving $\hat Q_j(\nu_l(k))$ are symmetric under a permutation of the indices of the different branches $\Omega_j$ of the dispersion relation, see Example~1 and 2, above. As a consequence, the integral over $\p D^+$ may be distributed to separate off these terms without introducing branching. This is to be expected: the solution of the original boundary-value problem should be independent of the choice of indices on $\{\Omega_j, j=1, \ldots, N\}$: switching our labels on $\Omega_1$ and $\Omega_2$, for instance, should not affect the solution. If the solution is to be symmetric under these permutations, we should expect the same for the terms containing $\hat Q_j(\nu_l(k))$.

\vs

{\bf Example 1.} We show that $R_1=0=R_2$ in \rf{kgfinal}. Recall that the original solution formulae \rf{kgsol} could not be deformed in the upper-half plane due to the presence of either $e^{-\Omega_1 t}$ or $e^{-\Omega_2 t}$. These exponentials canceled, and they are absent in \rf{r1r2}. On the other hand, the factor $e^{ikx}$ decays in the upper half plane and it follows from Jordan's Lemma that both $R_1$ and $R_2$ are zero. Thus \rf{kgfinal} with $R_1=0=R_2$ represent the final form of the solution of the KG equation posed on the positive half line with Dirichlet boundary data. The solution formula shows that no boundary information on $p(x,t)$ is required.

\vs

{\bf Example 2.} The fact that $S_1=0=S_2$ follows in exactly the same way: the region $D^+$ was previously inaccessible, but because the time-dependent exponential is absent, we may apply Cauchy's Theorem around $D^+$, with the integral contribution from the path at infinity vanishing, due to Jordan's Lemma. Thus $S_1=0=S_2$ since no singularities are present in $D^+$.

\vs\vs

{\bf Remarks.}

\begin{itemize}

\item For these particular example, one can observe that $R_1=0=R_2$, $S_1=0=S_2$ in another way. Consider $R_1$. Replacing $k\ra -k$ in the integral, we get

\begin{align}\nonumber
R_1&=\frac{1}{2\pi}\int_{-\infty}^\infty e^{ikx}\hat q(-k,t)dk\\\nonumber
&=-\frac{1}{2\pi}\int_{\infty}^{-\infty} e^{-ikx}\hat q(k,t)dk\\\nonumber
&=\frac{1}{2\pi}\int_{-\infty}^\infty e^{ik(-x)}\hat q(k,t)dk,
\end{align}
which is the inverse Fourier transform of the Fourier transform of $q(x,t)$, but evaluated at $-x$. Since $x$ is positive, it follows that $R_1=0$. A similar argument works for $R_2$, demonstrating that $R_1=0=R_2$ without contour deformation. For $S_1$ and $S_2$, we deform back to the real line, after which the above argument can be repeated. This method works for these examples because the lone nontrivial symmetry $\nu_2(k)=-k$ indicates the presence of a mirror symmetry in \rf{kg} and (\ref{fn}-b): the substitution $x\ra -x$ leaves both equations invariant and one expects that the half-line boundary-value problems for these equations may be solved using the method of images.

\item In the first example, we eliminated $g_1^{(1}$ and $g_1^{(2)}$, the time transforms of the Neumann data. Similary, in the second example we eliminated $g_0^{(1}$ and $g_0^{(2)}$. One may solve the first new global relations for $g_0^{(1)}$ and the second one for $g_1^{(2)}$, for instance. When one does so, the terms depending on the Fourier transforms $\hat Q_j(\nu_l(k))$ of the solution evaluated at $-k$ are not symmetric under permutation of the indices, as we expect.
    Nonetheless, as the reader easily verifies for the KG equation with $g_0^{(1)}$ and $g_0^{(2)}$ eliminated ({\em i.e.}, $g_0^{(2)}$ and $g_0^{(1)}$ are assumed known), the contributions from the term containing $\hat q(-k, t)$ still vanish, due to Jordan's Lemma. Thus the resulting solution formula for $q(x,t)$ and $p(x,t)$ is not symmetric under a permutation of the dispersion branch indices. This may be explained by noting that specifying boundary conditions by supplying $g_0^{(2)}$ and $g_0^{(1)}$ breaks the symmetry at the level of the problem statement already: when one specifies $g_0^{(2)}=g_0(\Omega_2,t)=\int_0^t e^{\Omega_2 s}q(0,s)ds$, one needs to include the information whether this function of $t$ originates from integration involving $\Omega_1$ or $\Omega_2$.

\end{itemize}

\section{A three-dimensional example due to Fokas and Treharne}

In \cite{treharnefokas}, Treharne and Fokas study a three-dimensional linear system on the half line, originating from a problem in elasticity. We discuss how their systems fits within the algorithm outlined in this paper. Although \cite{treharnefokas} is a seminal paper, as the first work where the UTM is considered for a system, the arguments presented there are specific to that system (the same can be said for the considerations in \cite{fokaspellonisys}) and do not highlight generic features of the UTM as applied to systems. In particular, we have emphasized the role played by expressions involving symmetric functions of the roots of the dispersion relations, allowing us to bypass the need for any branch cuts. Doing so, all expressions involved in the formal manipulations are analytic and much of the flavor of the scalar case is retained. Indeed the solution itself should be invariant under a relabeling of the branches of the dispersion relation, assuming the boundary conditions preserve this indifference.

The Treharne-Fokas system has the dispersion relation

\[
(i\lambda k^2+\omega)(\omega^2-k^2)-\alpha \beta \omega k^2=0.
\]

\no Here $\alpha$, $\beta$ and $\lambda$ are real parameters. This dispersion relation is cubic in $\omega$ and quartic in $k$ with a branch point at infinity in the $k$ variable. Other than by using Cardano's formulae, it is not possible to obtain the branches $\Omega_j(k)$ explicitly. The use of Cardano's formulae results, as expected, in horrendous expressions whose manipulation is not helpful.
However, by a straightforward asymptotic analysis of the dispersion relation, one can identify the regions in the complex $k$ plane for $\omega$ with simple asymptotic behavior as $k\rightarrow \infty$, allowing us to label the different branches $\Omega_1,\Omega_2,\Omega_3$. From the knowledge of the branches one readily obtains the associated eigenvectors (in terms of $\Omega_j$, $j=1,2,3$) that diagonalize the relevant operator. Much of the analysis presented above follows immediately and is not presented. Note that knowledge of the asymptotic behavior of $\Omega_j$, $j=1,2,3$ is sufficient to deform the paths of integration. At the next step one requires the symmetries that transmute one branch of the dispersion relation into another. Perhaps surprisingly, these symmetries are easily obtained in terms of $\omega$ from the dispersion relation: since the dispersion relation is a bi-quadratic in $k$, there are four symmetries: $\nu_1(k) = k$, $\nu_2(k) = -k$, $\nu_3(k) = f(\omega)$ ,$\nu_4(k) = -f(\omega)$,
where $f(\omega) = -\frac{1}{\sqrt \lambda} (\lambda \omega^2 + i(\alpha\beta + 1)\omega) $.  With the symmetries in hand, the arguments presented earlier follows and a solution expression is obtained. Of course, one does not require every combination of $\omega_j$ and $\nu_j$ and a judicious choice should be made to eliminate as many boundary values as possible. It is here that the peculiarities of the PDE problem at hand come into play.

The solution formula obtained as described above depends on the expressions for $\Omega_1$, $\Omega_2$ and $\Omega_3$: in effect we have reduced the solution of a PDE boundary-value problem to that of solving a third-order polynomial. For a general system, the solution expression depends on the roots $\Omega_j(k)$ and on the symmetries $\nu_j(k)$. Neither one can be written down explicitly for general dispersion relations. The complexity of solving boundary-value problems for linear systems of constant-coefficient PDEs is thus reduced to that of solving one additional polynomial problem compared to the scalar case. For scalar problems, $\Omega(k)$ is given explicitly, but the determination of the symmetries requires the solution of a polynomial of degree one less than the order of the problem (after the trivial symmetry $\nu_1(k)=k$ has been divided out). For systems, one faces the additional task of obtaining the branches of the dispersion relation $\Omega_j(k)$.

The analysis in \cite{treharnefokas} differs from ours in that the authors of \cite{treharnefokas} parameterize the dispersion relation, thereby diagonalizing the operator and eliminating any need to analyze expressions which may be branched. In our work, by employing symmetric functions of the branches of the dispersion relation, one yet again avoids the detailed analysis of branched expressions. The benefit of parametrization is that expressions are obtained that are more explicit. Of course, the fact that this problem may be suitably parametrized is not generic. Only dispersion relations given by an algebraic curve of genus zero can be parameterized. The genus of the dispersion relation and its corresponding parametrization may be readily obtained using the symbolic algebra package Maple and its package {\tt algcurves}.

\section{Conclusions}

We have demonstrated how Fokas's Unified Transform Method can be generalized from scalar evolution PDEs with constant coefficients to systems of such equations. Our main goal has been to show that the method continues its applicability much as in the scalar case, with added complexity relative to the system under consideration. It is difficult to make general statements for such a large class of problems, but some trends are clear. Although the different branches of the dispersion relation typically contain radicals, the different stages of the solution process can be executed without the need for branch cuts to be introduced, as the functions that arise are symmetric functions of the dispersion relation branches. Further, as in the scalar case, the symmetries of the dispersion relation give rise to new global relations, allowing for the elimination of unknown boundary data. Alternatively, this process allows one to determine the amount of boundary information that needs to be supplied to have a well-posed problem.

\section*{Acknowledgements}
Thomas Trogdon is acknowledged for useful conversations. This work was supported by the National
Science Foundation through grant NSF-DMS-1008001 (BD). This material is partially based upon work supported by the National Science Foundation under Grant No. DMS-1439786 while BD and VV were in residence at the Institute for Computational and Experimental Research in Mathematics in Providence, RI, during the Spring 2017 semester. Any opinions, findings,
and conclusions or recommendations expressed in this material are
those of the authors and do not necessarily reflect the views of the
funding sources.

\appendix

\section*{Appendix}

In this appendix, we show how to apply the UTM to the wave equation and equations like it, posed on the half line $x>0$. It is surprising that of the different PDEs typically dealt with in a first PDE course, the treatment of the wave equation using the UTM is not found in the literature. The closest to is the solution of the wave equation in a moving domain in \cite{pp1}, but the approach there requires knowledge of Lax pairs and Riemann-Hilbert problems. We show below that the use of the UTM leads to solution formulae that are valid in the whole quarter plane $x>0$, $t>0$. This is in contrast to the use of d'Alembert's formula, which leads to different solution representations in different regions of the quarter plane. Fourier transform methods may be used as well, but those methods run into difficulties if the boundary conditions are not homogeneous or if Robin boundary conditions are specified, for instance. We do not claim that the use of the UTM is the most convenient way to solve these wave equation-like problems, but it is instructional to see how the solution of such equations fits within the UTM framework.

We discuss these wave equations in an appendix rather than the body of the paper as the two solutions of the dispersion relation are not branched, and thus their solution method is not the typical application of the UTM outlined in the main sections. We examine the problem on $x>0$ for the wave-like equation

\beq\la{mm}
u_{tt}-a u_{xt}-u_{xx}=0.
\eeq

\no Here $a\in \mathbb{R}$. Note that the wave equation is recovered if $a=0$. We follow the algorithm for the UTM for systems and we rewrite this equation as

\alpheqn
\bea\la{mmsys}
u_{t}&=&v,\\
v_{t}&=&a v_{x}+u_{xx}.
\eea
\resetalpheqn

\no This leads to
\beq
\Lambda(k)=\twomatrix{0}{-1}{k^2}{-i a k},
\eeq
and
\beq
\Omega_{1,2}=i k \alpha_{1,2}=ik \frac{-a\pm \sqrt{4+a^2}}{2}.
\eeq
Note that $\alpha_1>0$, $\alpha_2<0$. For the wave equation, $\alpha_{1,2}=\pm 1$. Since this implies that one of the characteristic speeds is positive (information is carried away from the boundary), and one is negative (information is carried towards the boundary), we expect that one boundary condition needs to be prescribed.

With
\beq
X=\twomatrix{0}{0}{ik+\p_x}{a},
\eeq

\no the local relations are

\beq
\left(e^{-ikx+\Omega_j t}\left((\Omega_j+iak)u-v\right)\right)_t
+\left(e^{-ikx+\Omega_j t}\left(ik u+u_x+av\right)\right)_x=0, ~~~j=1,2.
\eeq

\no This leads to the global relations

\beq\la{mmgr}
(\Omega_j+iak)\hat u_0(k)-\hat v_0(k)-e^{\Omega_j t}(\Omega_j+iak)\hat u(k,t)+e^{\Omega_j t}\hat v(k,t)+a h_0^{(j)}+g_1^{(j)}+ik g_0^{(j)}=0,
\eeq
for $j=1,2$, valid in Im$\,k\leq 0$. Here

\begin{align}\nonumber
&&\hat u_0(k)&=\int_0^\infty e^{-ikx} u(x,0)dx, && \hat v_0(k)=\int_0^\infty e^{-ikx} v(x,0)dx,\\
&&\hat u(k,t)&=\int_0^\infty e^{-ikx} u(x,t)dx, && \hat v(k,t)=\int_0^\infty e^{-ikx} v(x,t)dx,\\\nonumber
&&g_0^{(j)}(t)&=\int_0^t e^{\Omega_j s}u(0,s)ds, && g_1^{(j)}(t)=\int_0^t e^{\Omega_j s}u_x(0,s)ds,\\\nonumber
&&h_0^{(j)}&=\int_0^t e^{\Omega_j s} v(0,s)ds=\int_0^t e^{\Omega_j s} u_t(0,s)ds.
\end{align}
It should be noted that specifying the Dirichlet boundary condition $u(0,t)$ results in both $g_0^{(j)}$ and $h_0^{(j)}$ being known, since $u_t(0,t)$ is obtained by taking a time derivative of $u(0,t)$.

From the Global Relations \rf{mmgr}, we obtain the solution formula

\beq\la{mmsol}
u(x,t)=I_1(x,t)+I_2(x,t),
\eeq

\no with

\alpheqn
\begin{align}
I_1(x,t)&=\frac{1}{2\pi}\int_{-\infty}^\infty e^{ikx}\frac{e^{-\Omega_2 t}\left((\Omega_2+iak)\hat u_0(k)-
\hat v_0(k)\right)-e^{-\Omega_1 t}\left((\Omega_1+iak)\hat u_0(k)-
\hat v_0(k)\right)}{\Omega_2-\Omega_1}dk,\\\la{i2}
I_2(x,t)&=\frac{1}{2\pi}\int_{-\infty}^\infty e^{ikx}\frac{e^{-\Omega_2 t}\left(a h_0^{(2)}+g_1^{(2)}+ikg_0^{(2)}\right)-e^{-\Omega_1 t}\left(a h_0^{(1)}+g_1^{(1)}+ikg_0^{(1)}\right)}{\Omega_2-\Omega_1}dk.
\end{align}
\resetalpheqn

\no We do not write down a solution formula for $v(x,t)=u_t(x,t)$. The first term in \rf{mmsol} is known and no manipulation of it is necessary. The second term brings in dependence from both the Dirichlet and Neumann boundary data and we expect to be able to eliminate one of these.

We wish to deform the path of integration used for $I_2$ as far as possible from the real line in the $k$ plane. As for the KG equation, $D^+=\{k\in \mathbb{C}: \mbox{Im}\,k>0\}$, thus no deformation away from the real line is allowed.

Next, we examine the symmetries of the dispersion relation. We solve

\beq
\det(\Lambda(\nu(k))-\omega(k)I)=0,
\eeq
where $\omega=\Omega_1$ or $\omega=\Omega_2$. We obtain

\beq
\nu_1(k)=k ~~(2\times),~~ \nu_2(k)=\frac{\alpha_1}{\alpha_2}k,~~ \nu_3(k)=\frac{\alpha_2}{\alpha_1}k.
\eeq

\no These symmetries are the reason we proceed with Example \rf{mm} instead of the regular wave equation: these symmetries are reminiscent of those used to solve scalar interface problems, see \cite{deconinckpellonisheils}, for instance. More importantly, they go beyond $k\ra -k$, which is used in the main sections of the paper.

Since $\Omega_1(\nu_3(k))=\Omega_2(k)$ and $\Omega_2(\nu_2(k))=\Omega_1(k)$, we obtain two new global relations, both valid on the real line: the first one is obtained by substituting $k\ra \nu_3(k)$ in \rf{mmgr} with $j=1$, while the second is obtained by letting $k\ra \nu_2(k)$ in \rf{mmgr} with $j=2$. These new global relations are

\newcommand{\aot}{\frac{\alpha_1}{\alpha_2}}
\newcommand{\ato}{\frac{\alpha_2}{\alpha_1}}

\alpheqn
\begin{align}\nonumber
\left(\Omega_2+ia\ato k\right)\hat u_0\left(\ato k\right)-\hat v_0\left(\ato k\right)
-e^{\Omega_2 t}\left(\Omega_2+ia\ato k\right)\hat u\left(\ato k,t\right)+\\\la{ngr1}
+e^{\Omega_2 t}\hat v\left(\ato k,t\right)+a h_0^{(2)}+g_1^{(2)}+i\ato k g_0^{(2)}=0,\\\nonumber
\left(\Omega_1+ia\aot k\right)\hat u_0\left(\aot k\right)-\hat v_0\left(\aot k\right)
-e^{\Omega_1 t}\left(\Omega_1+ia\aot k\right)\hat u\left(\aot k,t\right)+\\\la{ngr2}
+e^{\Omega_1 t}\hat v\left(\aot k,t\right)+a h_0^{(1)}+g_1^{(1)}+i\aot k g_0^{(1)}=0,
\end{align}
\resetalpheqn

\no and both are valid for Im\,$k\geq 0$, since $\alpha_1/\alpha_2<0$.

Next, we wish to eliminate either the Dirichlet or Neumann data from \rf{i2}. If Dirichlet data is given, we wish to use (\ref{ngr1}-b) to eliminate the Neumann data functions $g_1^{(j)}$, $j=1, 2$. On the other hand, if Neumann boundary conditions are provided, we use (\ref{ngr1}-b) to eliminate $g_0^{(j)}$, $j=1, 2$. As in the scalar case \cite{dtv, fokasbook}, if a Robin boundary condition is given, all of $g_m^{(j)}$, $j=1, 2$, $m=0, 1$ are eliminated using (\ref{ngr1}-b) and the two time transforms of the Robin condition. For exposition sake, let us assume Dirichlet data is given. We substitute $g_1^{(j)}$, $j=1, 2$ obtained from (\ref{ngr1}-b) in \rf{i2}. This results in

\beq
I_2(x,t)=T_1(x,t)+T_2(x,t)+T_3(x,t),
\eeq
where

\alpheqn
\begin{align}\nonumber
T_1(x,t)=&\frac{1}{2\pi}\int_{-\infty}^\infty \frac{e^{ikx}}{\Omega_2-\Omega_1}
\left(
e^{-\Omega_2 t}\left(\hat v_0\left(\ato k\right)-\left(\Omega_2+ia\ato k\right)\hat u_0\left(\ato k\right)\right)+\right.\\
&\left.-e^{-\Omega_1 t}\left(\hat v_0\left(\aot k\right)-\left(\Omega_1+ia\aot k\right)\hat u_0\left(\aot k\right)\right)
\right)dk, \\
T_2(x,t)=&\frac{1}{2\pi}\int_{-\infty}^\infty \frac{ik e^{ikx}}{\Omega_2-\Omega_1}\left(
e^{-\Omega_2 t} \left(1-\ato\right)g_0^{(2)}-e^{-\Omega_1 t} \left(1-\aot\right)g_0^{(1)}
\right)dk,\\\nonumber
T_3(x,t)=&\frac{1}{2\pi}\int_{-\infty}^\infty \frac{e^{ikx}}{\Omega_2-\Omega_1}\left(
\left(\Omega_2+ia\ato k\right)\hat u\left(\ato k, t\right)-\hat v\left(\ato k, t\right)+\right.\\
&\left.-\left(\Omega_2+ia\ato k\right)\hat u\left(\ato k, t\right)+\hat v\left(\ato k, t\right)
\right)dk.
\end{align}
\resetalpheqn

\no Note that $T_1(x,t)$ is determined by the initial data. It can be combined with $I_1(x,t)$. The function $T_2(x,t)$ depends only on the Dirichlet data, assumed known. Lastly, $T_3(x,t)$ exhibits dependence on the solution, but we show below that $T_3(x,t)\equiv 0$. Thus

\beq\la{mmfinal}
u(x,t)=I_1(x,t)+T_1(x,t)+T_2(x,t),
\eeq

\no is a solution of \rf{mm}, posed on $x>0$ with Dirichlet data specified at $x=0$. This solution is valid for all $x$ and $t$ in the quarter plane $x>0$, $t>0$.

To show that $T_3(x,t)\equiv 0$, it suffices to note that due to the absence of the exponential, deformation of the path of integration into $D^+$ (the upper-half plane) is allowed, and the conclusion follows from Jordan's Lemma and Cauchy's Theorem. Observe that the singularity at $k=0$ is removable.

\begin{figure}[tb]
\begin{center}
\def\svgwidth{5in}
    \vspace*{0in}
    \hspace*{0.0in}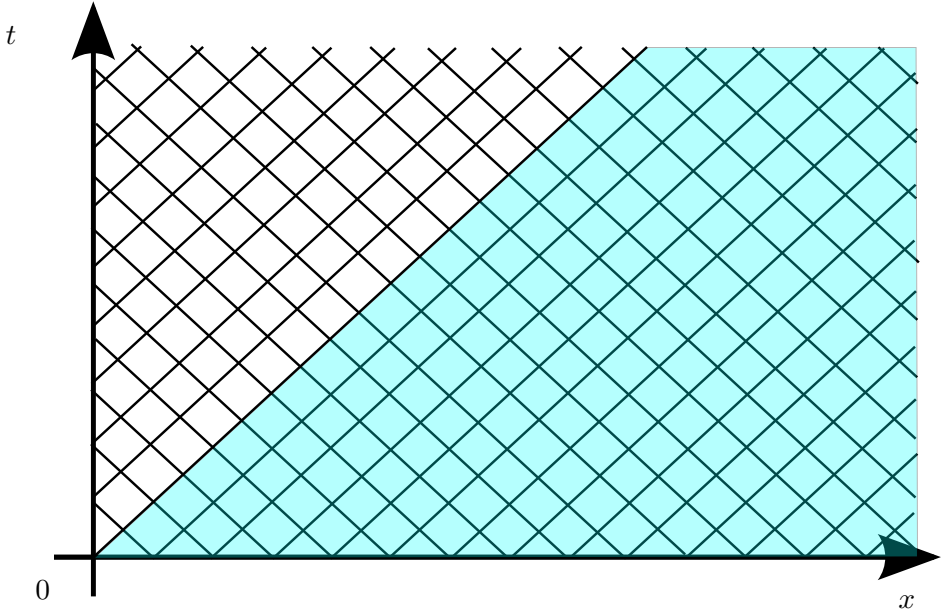
    \vspace*{0in}
    \caption{\la{fig:wavechar} In the blue region, the solution depends only on initial data. Above this region, the solution depends on both initial and boundary data.}
\end{center}
\end{figure}

Lastly, we restrict to the case of the wave equation ($a=0$, $\Omega_{1,2}=\pm ik$) and we illustrate how the solution obtained here can be reduced to d'Alembert form. We examine the Dirichlet, Neumann, and Robin problem. The characteristics for the wave equation are drawn in Fig.~\ref{fig:wavechar}. One family of characteristics moves initial data away from the line $t=0$, to the left. The second family takes information to the right, either initial data coming from $t=0$ or else boundary data from $x=0$. It follows that points in the blue region in Fig.~\ref{fig:wavechar} depend on initial data only. Their solution expression should be the same, independent of which boundary conditions are imposed. Points above the blue region depend on both initial data and boundary data. In this region we expect a different solution expression, depending on which boundary conditions are specified.

Examining $I_1(x,t)$ with $a=0$, we find

\begin{align}\nonumber
I_1(x,t)&=\frac{1}{2\pi}\int_{-\infty}^\infty e^{ikx}\frac{e^{-\Omega_2 t}\left(\Omega_2\hat u_0(k)-
\hat v_0(k)\right)-e^{-\Omega_1 t}\left(\Omega_1\hat u_0(k)-
\hat v_0(k)\right)}{\Omega_2-\Omega_1}dk\\\nonumber
&=\frac{1}{2\pi}\int_{-\infty}^\infty e^{ikx} \frac{-ik \hat u_0(k)\left(e^{ikt}-e^{-ikt}\right)-\hat v_0(k)\left(e^{ikt}-e^{-ikt}\right)}{-2ik}dk\\\nonumber
&=\frac{1}{4\pi}\int_{-\infty}^\infty \hat u_0(k)\left(
e^{ik(x+t)-e^{ik(x-t)}}
\right)dk+\frac{1}{4i\pi}\int_{-\infty}^\infty \hat v_0(k)\frac{e^{ik(x+t)}-e^{ik(x-t)}}{k}dk\\\la{dalembert}
&=\frac{1}{2}\left(u_0(x+t)+u_0(x-t)\right)+\frac{1}{2}\int_{x-t}^{x+t}v_0(s)ds,
\end{align}
using well-known calculational properties of the Fourier transform. This is the d'Alembert form of the solution, as expected in the blue region $x>t$. This implies that we should find $I_2(x,t)=0$ in the blue region. Indeed, with $a=0$ and $x>t$,

\begin{align}\nonumber
I_2(x,t)&=\frac{1}{2\pi}\int_{-\infty}^\infty e^{ikx}\frac{e^{-\Omega_1 t}\left(g_1^{(1)}+ikg_0^{(1)}\right)-e^{-\Omega_2 t}\left(g_1^{(2)}+ikg_0^{(2)}\right)}{\Omega_1-\Omega_2}dk\\\nonumber
&=\frac{1}{4ik}\int_{-\infty}^\infty e^{ikx}\frac{e^{-ikt}\left(g_1^{(1)}+ikg_0^{(1)}\right)-e^{ikt}\left(g_1^{(2)}+ikg_0^{(2)}\right)}{k}dk\\
&=\frac{1}{4\pi i}\int_{-\infty}^\infty dk \int_0^t ds \frac{\left(e^{ik(x-t+s)}-e^{ik(x+t+s)}\right)(u_x(0,s)+ik u(0,s))}{k}.
\end{align}

\no Both $x-t+s>0$ and $x+t+s>0$ if $x>t$. Thus the integrands are decaying exponentially in the upper-half $k$ plane and by Cauchy's Theorem and Jordan's Lemma, $I_2(x,t)\equiv 0$. Note that the singularity at $k=0$ in the above integrals is removable. Thus in $x>t$,

\beq\la{dalembert1}
u(x,t)=\frac{1}{2}\left(u_0(x+t)+u_0(x-t)\right)+\frac{1}{2}\int_{x-t}^{x+t}v_0(s)ds.
\eeq

\no This formula may be extended to include $x=t$, provided the initial and boundary data are compatible at $(x,t)=(0,0)$. Next, we examine $x<t$.

For the {\bf Dirichlet problem}, $I_2(x,t)=T_1(x,t)+T_2(x,t)$. With $a=0$,

\begin{align}\nonumber
T_1(x,t)&=\frac{1}{2\pi}\int_{-\infty}^\infty \frac{e^{ikx}}{\Omega_1-\Omega_2}
\left(
e^{-\Omega_1 t}\left(\hat v_0\left(-k\right)-\Omega_1\hat u_0\left(-k\right)\right)-e^{-\Omega_2 t}\left(\hat v_0\left(-k\right)-\Omega_2\hat u_0\left(-k\right)\right)
\right)dk\\\nonumber
&=\frac{1}{2\pi}\int_{-\infty}^\infty \frac{
\hat v_0(-k)\left(e^{ik(x-t)}-e^{ik(x+t)}\right)-ik\hat u_0(-k)\left(e^{ik(x-t)}-e^{ik(x+t)}\right)
}{2ik}dk\\
&=-\frac{1}{2} u_0(t-x)-\frac{1}{2}\int_0^{t-x}v_0(s)ds.
\end{align}

Next, again with $a=0$,
\begin{align}\nonumber
T_2(x,t)&=\frac{1}{\pi}\int_{-\infty}^\infty \frac{ik e^{ikx}}{\Omega_1-\Omega_2}\left(
e^{-\Omega_1 t} g_0^{(1)}-e^{-\Omega_2 t} g_0^{(2)}
\right)dk\\\nonumber
&=\frac{1}{2\pi}\int_{-\infty}^\infty dk \int_0^t \!ds~ u(0,s)\left(e^{ik(x-t+s)}-e^{ik(x+t-s)}\right)\\
&=u(0,t-x),
\end{align}
where we have used the definition of the delta function and the fact that the Dirichlet boundary condition $u(0,t)=0$ for $t<0$. Using that $u_0(x)=0$, $v_0(x)=0$ for $x<0$ in \rf{dalembert}, we finally obtain

\beq
u(x,t)=\left\{
\ba{lcr}
\ds \frac{1}{2}(u_0(x+t)-u_0(t-x))+\frac{1}{2}\int_{t-x}^{x+t}v_0(s)ds+u(0,t-x),&x<t,\\
\ds \frac{1}{2}(u_0(x+t)+u_0(x-t))+\frac{1}{2}\int_{x-t}^{x+t}v_0(s)ds, & x>t.
\ea
\right.
\eeq

\no It is a straightforward check that this solution satisfies the wave equation, the initial conditions, and the Dirichlet boundary condition. Further, if initial and boundary data are compatible at $(x,t)=(0,0)$ then the solution is continuous at $x=t$.

For the {\bf Neumann problem}, $I_2(x,t)=T_1(x,t)+T_2(x,t)$. With $a=0$ and $x<t$, the new global relations (\ref{ngr1}-b) are

\alpheqn
\begin{align}
ik g_0^{(1)}&=ik \hat u_0(-k)-\hat v_0(-k)-ik e^{ikt}\hat u(-k,t)+e^{ikt}\hat v(-k,t)+g_1^{(1)},\\
ik g_0^{(2)}&=-ik \hat u_0(-k)-\hat v_0(-k)+ik e^{-ikt}\hat u(-k,t)+e^{-ikt}\hat v(-k,t)+g_1^{(2)},
\end{align}
\resetalpheqn
which are substituted in \rf{i2}, resulting in

\begin{align}\nonumber
I_2(x,t)&=\frac{1}{2\pi}\int_{-\infty}^\infty e^{ikx} \frac{e^{ikt}\left(g_1^{(2)}+ik g_0^{(2)}\right)-e^{-ikt}\left(g_1^{(1)}+ik g_0^{(1)}\right)}{-2ik}dk\\\nonumber
&=-\frac{1}{2\pi}\int_{-\infty}^\infty \frac{e^{ik(x+t)}g_1^{(2)}-e^{ik(x-t)}g_1^{(1)}}{ik}dk+\\\nonumber
&~~~-\frac{1}{4\pi} \int_{-\infty}^\infty
\frac{-ik \hat u_0(-k) \left(e^{ik(x+t)}+e^{ik(x-t)}\right)+\hat v_0(-k) \left(e^{ik(x-t)}-e^{ik(x+t)}\right)}{ik} dk\\\nonumber
&=-\frac{1}{2\pi}\int_{-\infty}^\infty dk \int_0^t ds\,\, u_x(0,s)\frac{e^{ik(x+t-s)}-e^{ik(x-t+s)}}{ik}+\\\nonumber
&~~~-\frac{1}{4\pi} \int_{-\infty}^\infty
\frac{ik \hat u_0(k) e^{ik(t-x)}+\hat v_0(k) \left(e^{ik(t-x)}-e^{-ik(x+t)}\right)}{-ik} dk\\
&=-\int_0^t u_x(0,s)\left(\theta(x+t-s)-\theta(x-t+s)\right)ds+\frac{1}{2}u_0(t-x)+\frac{1}{2}\int_{0}^{t-x} v_0(s)ds
\end{align}

\no We have omitted the terms containing $\hat u(-k,t)$ and $\hat v(-k, t)$, which are easily shown to have zero contributions. Here

\beq
\theta(x)=\left\{
\ba{lcr}
\ds 0, &x<0,\\
\ds 1, & x>0,
\ea
\right.
\eeq

\no the Heaviside function. As before, we have used that $u_0(x)=0$, $v_0(x)=0$ for $x<0$. Combining our results, we obtain

\beq
u(x,t)=\left\{
\ba{lcr}
\ds \frac{1}{2}(u_0(x+t)+u_0(t-x))+\frac{1}{2}\left(\int_0^{t-x}v_0(s)ds+\int_0^{x+t}v_0(s)ds\right)+\\
\ds ~~~~~~~~~~~~~~~~~~~~~~~~~~~~~~~~~~~~~~~~~~~~~~~~-\int_0^{t-x} u_x(0,s)ds,&x<t,\\
\ds \frac{1}{2}(u_0(x+t)+u_0(x-t))+\frac{1}{2}\int_{x-t}^{x+t}v_0(s)ds, & x>t.
\ea
\right.
\eeq

\no As for the Dirichlet problem, it is straightforward to check that this solution satisfies the wave equation, the initial conditions, and the Neumann boundary condition. If initial and boundary data are compatible at $(x,t)=(0,0)$ then the solution is continuous at $x=t$.

Lastly, we consider {\bf Robin boundary conditions:}
\beq
a u(0,t)+b u_x(0,t)=f(t),
\eeq

\no with given real constants $a$ and $b$ and a time-dependent function $f(t)$. We evaluate this boundary condition at $t=s$, multiply by $\exp{\Omega_j s}$ and integrate, to obtain

\beq\la{robintt}
\gamma g_0^{(j)}+g_1^{(j)}=f^{(j)}, ~~~~j=1, 2,
\eeq

\no with

\beq
f^{(j)}=f(\Omega_j, t)=\int_0^t e^{\Omega_j s} f(s)ds, ~~~j=1,2.
\eeq

\no The two equations \rf{robintt} are valid for all $k\in \mathbb{C}$. We solve these equations combined with the new global relations (\ref{ngr}a-b) for $g_m^{(j)}$, $j=1, 2$, $m=0, 1$, obtaining

\alpheqn
\begin{align}\la{robing}
g_0^{(j)}&=\frac{f^{(j)}+(\Omega_j \hat u_0(-k)-\hat v_0(-k))-e^{\Omega_j t}(\Omega_j \hat u(-k,t)-\hat v(-k,t))}{\gamma+ik},\\
g_1^{(j)}&=\frac{ik f^{(j)}-\gamma (\Omega_j \hat u_0(-k)-\hat v_0(-k))+\gamma e^{\Omega_j t} (\Omega_j \hat u(-k,t)-\hat v(-k,t))}{\gamma+ik}.
\end{align}
\resetalpheqn

\no These expressions are substituted in \rf{i2} with $a=0$. This results in

\begin{align}\la{i2robin}
I_2(x,t)=J_1(x,t)+J_2(x,t)+J_3(x,t),
\end{align}

\no with

\beq
J_1(x,t)=\frac{1}{2\pi}\int_C \frac{e^{ikx}}{\gamma+ik}\left(f^{(1)}e^{-ikt}-f^{(2)}e^{ikt}\right)dk,
\eeq

\no which depends on the boundary conditions only,

\beq
J_2(x,t)=-\frac{1}{4\pi}\int_C \frac{e^{ikx}(\gamma-ik)}{\gamma+ik}\hat u_0(-k)\left(e^{-ikt}+e^{ikt}\right)dk,
\eeq

\no the dependence on the initial condition $u_0(x)$, and

\beq
J_3(x,t)=\frac{1}{4\pi}\int_C \frac{e^{ikx}(\gamma-ik)}{ik(\gamma+ik)}\hat v_0(-k)\left(e^{-ikt}-e^{ikt}\right)dk,
\eeq

\no the contribution from the initial condition $v_0(x)$. In these definitions, the contour $C$ accommodates the presence of the singularity at $k=i\gamma$: in the integral obtained after substituting (\ref{robing}-b) in \rf{i2}, the singularity at $k=i\gamma$ is removable, since the integrand of \rf{i2} is analytic in all of $\mathbb{C}$. If we distribute the integral so as to isolate the dependence of $\hat u(-k,t)$ (the term with $\hat v(-k, t)$ cancels), both integrals have a pole singularity. If $\gamma>0$, this singularity is in the upper-half plane, which is where we wish to close the contour to show that the contribution from $\hat u(-k, t)$ vanishes. To avoid this issue, we deform the path of integration to pass above $k=i\gamma$ along  path $C$, prior to distributing the integral. Thus $C$ is the real line if $\gamma<0$, and $C$ is asymptotic to the real line but passes above $i\gamma$ if $\gamma>0$. Using this contour, it immediately follows that the terms containing $\hat u(-k, t)$ vanish, as desired. Thus we have omitted these terms in \rf{i2robin}.

In what follows, we need

\begin{align}\la{ci1}
\frac{1}{2\pi}\int_{-\infty}^\infty \frac{1}{\gamma+ik} e^{ikx}dk&=\sgn(x)\theta(\gamma x)e^{-\gamma x},\\\la{ci2}
\frac{1}{2\pi}\int_{-\infty}^\infty \frac{\gamma-ik}{\gamma+ik}e^{ikx}dk&=-\delta(x)+2\,\sgn(x)\gamma \theta(\gamma x) e^{-\gamma x},
\end{align}

\no which are obtained using the definition of the delta function and standard contour integration. As elsewhere in this paper, all integrals over $k$ are principal-value integrals.  Using the definition of $f^{(j)}$, switching the order of integration, and using \rf{ci1}, we find

\alpheqn
\begin{align}
J_1(x,t)&=-\int_0^{t-x}f(s)e^{\gamma(t-x-s)}ds, \\
J_2(x,t)&=\frac{1}{2}u_0(t-x)+\gamma \int_0^{t-x} u_0(y)e^{\gamma(t-x-y)}dy,\\
J_3(x,t)&=-\frac{1}{2}\int_0^{t-x}v_0(y)dy+\int_0^{t-x} v_0(y)e^{\gamma (t-x-y)}dy,
\end{align}
\resetalpheqn

\no for both $\gamma>0$ and $\gamma<0$. For $\gamma>0$, these results include the integral contribution from the real line and a subtracted residue contribution as a consequence of deforming $C$ back to the real line.

Combining our results, we obtain

\beq
u(x,t)=\left\{
\ba{lcr}
\ds \frac{1}{2}(u_0(x+t)+u_0(t-x))+\gamma \int_0^{t-x} u_0(y)e^{\gamma(t-x-y)}dy+\\
~~\ds +\frac{1}{2}\int_{t-x}^{t+x}v_0(y)dy+\int_0^{t-x}v_0(s)e^{\gamma(t-x-y)}dy-\int_0^{t-x} f(s)e^{\gamma(t-x-s)}ds,&x<t,\\
\ds \frac{1}{2}(u_0(x+t)+u_0(x-t))+\frac{1}{2}\int_{x-t}^{x+t}v_0(s)ds, & x>t,
\ea
\right.
\eeq

\no the d'Alembert form of the solution for the half-line boundary-value problem for the wave equation with Robin boundary data. A direct calculation verifies that $u(x,t)$ satisfies the Robin boundary condition, and is continuous along $x=t$.

\bibliographystyle{plain}

\bibliography{mybib}

\end{document}